\documentclass[11pt]{article}
\usepackage{amsmath,amsxtra}
\usepackage{mathtools}
\usepackage{amssymb}
\usepackage{mathrsfs}
\usepackage{amsfonts}
\usepackage{esint}
\usepackage{amscd}
\usepackage[dvipdf]{graphics}
\usepackage[dvips]{graphicx}
\usepackage{graphicx,subfigure}
\usepackage{latexsym}
\usepackage{epstopdf}
\usepackage{latexsym}
\usepackage{amsmath}
    \def\dfrac{\displaystyle\frac}
    
\newcommand{\beqa}{\begin{eqnarray}}
\newcommand{\eeqa}{\end{eqnarray}}
\newcommand{\ba}{\begin{eqnarray*}}
\newcommand{\ea}{\end{eqnarray*}}

\date{}

\newtheorem{definition}{Definition}

\def\rw{\rightarrow}
\def\n{{\bf n}}
\def\var{\varphi}
\def\div{{\,\rm div \,}}

\def\tan{{\,\rm tan \,}}
\def\id{{\,\rm id \,}}

\def\so{{\,\rm so \,}}
\def\n{\mathbf{n}}
\def\sym{{\,\rm sym \,}}
\def\dist{{\,\rm dist \,}}
\def\id{{\,\rm id \,}}

\def\lam{\lambda}

\def\R{\mathbb{R}}

\def\skew{{\,\rm skew \,}}
\def\SO{{\,\rm SO \,}}

\def\lam{\lambda}
\def\be{{\begin{equation}}}
\def\ee{{\end{equation}}}

\def\Ga{{\Gamma}}
\def\pl{{\partial}}
\def\qfq{{\quad\mbox{for}\quad}}
\def\beq{\arraycolsep=1.5pt\begin{eqnarray}}
\def\eeq{\end{eqnarray}}

\newfont{\Blackboard}{msbm10 scaled 1200}

\newfont{\roma}{cmr10 scaled 1200}

\def\qfq{{\quad\mbox{for}\quad}}
\def\<{{\langle}}
\def\>{{\rangle}}

\newtheorem{thm}{{}\hskip\parindent Theorem}[section]
\newtheorem{lem}{{}\hskip\parindent Lemma}[section]
\newtheorem{pro}{{}\hskip\parindent Proposition}[section]

\newtheorem{rem}{{}\hskip\parindent Remark}[section]

\def\dfrac{\displaystyle\frac}

\def\be{\begin{equation}}
\def\ee{\end{equation}}
\def\beq{\arraycolsep=1.5pt\begin{eqnarray}}
\def\eeq{\end{eqnarray}}

\large
\title{\bf The time-dependent von K\'{a}rm\'{a}n shell equation as a limit of three-dimensional nonlinear elasticity }
\date{}

\author{Yizhao Qin\quad
Peng-Fei Yao\thanks{Corresponding author.\ Email: pfyao@iss.ac.cn}\quad \\[0.3cm]
Key Laboratory of  Systems and Control\\
Institute of Systems Science,
Academy of Mathematics and Systems Science\\
Chinese Academy of Sciences, Beijing 100190, P. R. China\\
School of Mathematical Sciences\\
University of Chinese Academy of Sciences, Beijing 100049, China
}

\begin{document}
\maketitle
\footnote{This work is  supported by the National
Science Foundation of China, grants  no. 61473126 and no. 61573342, and Key Research Program of Frontier Sciences, CAS, no. QYZDJ-SSW-SYS011.}
\begin{quote}
\begin{small}
{\bf Abstract} \,\,\, The asymptotic behaviour of solutions of three-dimensional nonlinear elastodynamics in a thin shell is considered, as the thickness $h$ of the shell tends to zero. Given the appropriate scalings of the applied force and of the initial data in terms of $h,$ it's verified that three-dimesional solutions of the nonlinear elastodynamic equations converge to solutions of the time-dependent von K\'{a}rm\'{a}n equations or dynamic linear equations for shell of arbitrary geometry.
\\[3mm]
\\[3mm]
\end{small}
\end{quote}

\section{Introduction and Main Results}
\setcounter{equation}{0}
\hskip\parindent In this paper, we concern about the rigorous derivation of the two-dimensional dynamic models for a thin elastic shell starting from three-dimensional nonlinear elastodynamics. To be clear, we consider a thin elastic shell of reference configuration
\begin{equation*}
S^h=\{z=x+s\mathbf{n}(x): x\in S,\quad -\frac{h}{2}<s<\frac{h}{2}\},\quad 0<h\leqslant h_0.
\end{equation*}
It's a family of shells of small thickness $h$ around the middle surface $S,$ where $S$ is a compact, connected, oriented 2d surface of the class $\mathcal{C}^2$ embedded in $\mathbb{R}^3$ with a $\mathcal{C}^2$ boundary $\partial S.$ By $\mathbf{n}(x),$ we denote the unit normal to $S$ and $S_x$ stands for the tangent space at $x.$ We suppose that
the energy potential of this thin shell $W:\mathbb{R}^{3\times3}\rightarrow[0,\infty]$ is a continuous function with the following properties:
\begin{eqnarray}
&&W(RF)=W(F),\quad \forall\quad F\in\mathbb{R}^{3\times3},\quad R\in \SO(3)(\mbox{frame indifference});\label{eq1.1} \\
&&W(R)=0,\quad \forall R\in \SO(3);\label{eq1.2} \\
&&\exists\quad\mbox{a positive constant}\quad C\quad\mbox{such that} \nonumber\\
&&\quad W(F)\geqslant C\dist^2(F,\SO(3)), \quad\forall F\in\mathbb{R}^{3\times3}; \label{eq1.3}\\
&&W\quad \mbox{is}\quad \mathcal{C}^2\quad\mbox{in a neighbourhood of}\quad \SO(3)\label{eq1.4};\\
&&\mid DW(F)\mid\leqslant C(\mid F\mid+1),\quad\forall F\in\mathbb{R}^{3\times 3}\label{eq1.23}.
\end{eqnarray}
Here, $\SO(3)$ denotes the group of proper rotations.
The dynamic equations of nonlinear elasticity arise from the action functional
\begin{equation*}
E^h(u^h)=\frac{1}{h}\int_0^{\xi_h}\int_{S^h}[\frac{\mid u^h_{\xi}\mid^2}{2}-W(\nabla u^h(z))+\langle f^h,u^h(z)\rangle]d\xi dz
\end{equation*}
and by computing the Euler-Lagrange equations of the above energy, the equations of elastodynamic read as
\begin{equation}\label{eq1.5}
\partial^2_{\xi}u^h-\div DW(\nabla u^h)=f^h\quad\mbox{in}\quad (0,\xi_h)\times S^h,
\end{equation}
where $u^h:[0,\xi_h]\times S^h\rightarrow\mathbb{R}^3$ is the deformation of the shell and $f^h:[0,\xi_h]\times S^h\rightarrow\mathbb{R}^3$ is an external body force applied to the shell. Equation (\ref{eq1.5}) is supplemented by the initial data
\begin{equation}\label{eq1.6}
u^h\mid_{\xi=0}=\bar{w}^h,\quad \partial_{\xi}u^h\mid_{\xi=0}=\hat{w}^h,
\end{equation}
and, respectively, by the mixed Neumann-clamped boundary conditions:
\begin{eqnarray}
&&u^h=z\quad\mbox{on}\quad (0,\xi_h)\times\{z=x+s\mathbf{n}(x): x\in\partial S,\, s\in(-\frac{h}{2},\frac{h}{2})\},\label{eq1.7} \\
&&DW(\nabla u^h)\mathbf{n}=0\quad\mbox{on}\quad (0,\xi_h)\times\{z=x\pm\frac{h}{2}\mathbf{n}(x): x\in S\}.\label{eq1.8}
\end{eqnarray}

Our purpose of this paper is to characterize the asymptotic behaviour of the solutions to (\ref{eq1.5}), as the thickness $h$ approaches to zero, by identifying the two-dimensional dynamic model for the thin elastic shell satisfied by their limit as $h\rightarrow 0.$

Lower dimensional models for thin bodies attract much attention in elasticity theory, as they are usually easier to handle from both analytical and numerical view than their three-dimensional counterparts. The problem of their rigorous derivation beginning from three-dimensional theory is one of the central issues in nonlinear elasticity.
In the stationary case, the application of variational methods, especially the $\Gamma-$ convergence, leads to the rigorous derivation of a hierarchy of limiting theories for thin plates and shells recently(\cite{FGM1}, \cite{FGM2}, \cite{LMP1}-\cite{LP}, \cite{yao}). The $\Gamma-$ convergence approach implies the convergence of minimizers of a sequence of functionals, to the minimizers of the limit. However, it doesn't guarantees the convergence of the possibly non-minimizing critical points (the equilibria), which are the solutions of the Euler-Lagrange equations of the corresponding functionals. In this setting, S. M\"{u}ller have first obtained convergence consequences in the von K\'{a}rm\'{a}n case for the thin plates in \cite{MP}. Then, the results of convergence of equilibria have been generalized to the cases of rods, beams and shells( \cite{MM1}, \cite{MMS} and \cite{L}). Under the physical growth condition of energy density, similar results are also established in \cite{MS}. For more detailed survey in this direction, see \cite{M}.

As for the time-dependent cases, the model from 3d to 2d has only been established when the energy per unit volume decays like
$h^4$ or stronger for thin plate so far(\cite{AMM1, AMM2}). Here we shall combine \cite{AMM2} and \cite{LMP2} to obtain the time-dependent model for the thin shells in the von K\'{a}rm\'{a}n case.

 Let $\Pi(x)=\nabla\mathbf{n}(x)$ denote the negative second fundamental form of $S$ at $x.$  Let $\pi$ be the projection onto $S$ along $\mathbf{n}(x),$ that is, $\pi(z)=x,$ for all $z=x+s\mathbf{n}(x)\in S^h.$ We assume that $0<h\leqslant h_0,$ with $h_0>0$ given sufficiently small to have $\pi$ well defined on each $S^h$ and $\frac{1}{2}<\mid \id+s\Pi(x)\mid<\frac{3}{2}$ for all $|s|<h_0/2.$

We recall some notations and results in the stationary case briefly. For a $H^1$ deformation $u,$ we associate its elastic energy (scaled per unit thickness) with
\begin{equation*}
I^h(u)=\frac{1}{h}\int_{S^h}W(\nabla u(z))dz,
\end{equation*}
where $W$ satisfy (\ref{eq1.1})-(\ref{eq1.4}) as well.
Furthermore, the total energy of thin shell in the stationary case is provided by
\begin{equation}\label{eq1.9}
J^h(u)=I^h(u)-\frac{1}{h}\int_{S^h}\langle f^h,u(z)\rangle dz,
\end{equation}
where the external force $f^h,$ defined on $S^h,$ is supposed to be
\begin{equation}\label{eq1.10}
f^h(x+s\mathbf{n})=h\sqrt{e^h}f(x)\det(\id+s\Pi(x))^{-1},\quad f(x)\in L^2(S,\mathbb{R}^3),\quad
\int_Sf(x)=0.
\end{equation}
 In (\ref{eq1.10}) $e^h>0$ is a given sequence obeying a prescribed scaling law. It's shown that if $f^h$ scales like $h^{\alpha},$ then the minimizers $u^h$ of $J^h(u)$ satisfy $I^h(u^h)\sim h^{\beta}$ with $\beta=\alpha$ if $0\leqslant\alpha\leqslant 2$ and $\beta=2\alpha-2$ if $\alpha>2.$ Throughout this note we shall assume that $\beta\geqslant 4,$ or more generally
\begin{equation}\label{eq1.11}
\lim_{h\rightarrow 0}\frac{e^h}{h^4}=\kappa<\infty.
\end{equation}
In particular, the case that $S\subset \mathbb{R}^2$ corresponding to the von K\'{a}rm\'{a}n and purely linear theories of plates is derived rigorously in \cite{FGM2}.

Let $\mathcal{V}(S,R^3)$  be the space of all $H^2$ infinitesimal isometries on $S.$ For each $V\in \mathcal{V}(S,R^3)$ for which there exists a matrix field $A\in H^1(S,\mathbb{R}^{3\times 3})$ such that
\begin{equation}\label{eq1.12}
\partial_{\tau}V(x)=A(x)\tau\quad\mbox{and}\quad A(x)^T=-A(x)\quad\forall x\in S,\quad \mbox{a.e.},\quad\tau\in S_x.
\end{equation}

For $F\in L^2(S,\mathbb{R}^{3\times 3}),$ let $F_{\tan}(x)=
[\langle F(x)\tau,\eta\rangle]_{\tau\eta\in S_x}.$
The quadratic forms $\mathcal{Q}_2(x,.)$ are given by
\begin{equation*}
\mathcal{Q}_2(x,F_{\tan})=\min\{\mathcal{Q}_3(\tilde{F}):\quad(\tilde{F}-F)_{\tan}=0\},\quad \mathcal{Q}_3(F)=D^2W(\id)(F,F).
\end{equation*}
The form $\mathcal{Q}_3$ is defined for all $F\in\mathbb{R}^{3\times 3},$ while $\mathcal{Q}_2(x,\cdot)$ for a given $x\in S,$ is defined on
tangential minors $F_{\tan}$ of $F\in\mathbb{R}^{3\times 3}.$ Both forms depend only on the symmetric parts of their arguments and are positive definite on the space of symmetric matrices(see \cite{FGM1}).

We define
the linear operators $\mathcal{L}_3:$ $\R^{3\times3}\rightarrow\R^{3\times3}$  and $\mathcal{L}_2(x,.):$ $\R^{2\times2}\rightarrow\R^{2\times2}$  by
\begin{equation*}
Q_3(F)=\mathcal{L}_3F:F\quad\mbox{and}\quad \mathcal{Q}_2(x,F_{\tan})=\mathcal{L}_2(x,F_{\tan}):F_{\tan}\quad\forall F\in\mathbb{R}^{3\times 3},
\end{equation*}
respectively,  $F_1:F_2=tr(F_1^TF_2)$ for two matrices $F_1$ and $F_2.$

If $\kappa=0$ in (\ref{eq1.11}), the $\Ga$-limiting of (\ref{eq1.9}) is given by
\begin{equation}\label{eq1.13}
J(V,\bar{Q})=\frac{1}{24}\int_S \mathcal{Q}_2(x,(\nabla(A\mathbf{n})-A\Pi)_{\tan})dx-\int_S\langle f,\bar{Q}V\rangle dx,\quad\forall V\in\mathcal{V},\quad\bar{Q}\in \SO(3).
\end{equation}
For $\kappa>0,$ the $\Gamma-$limit  of $J^h$ is
\beq\label{eq1.14}
J(V,B_{\tan},\bar{Q})&&=\frac{1}{2}\int_S \mathcal{Q}_2(x,B_{\tan}-\frac{\sqrt{\kappa}}{2}(A^2)_{\tan})dx+\frac{1}{24}\int_S \mathcal{Q}_2(x,(\nabla(A\mathbf{n})-A\Pi)_{\tan})dx \nonumber\\
&&-\int_S\langle f,\bar{Q}V\rangle dx,
\eeq
where $B_{\tan}$ on $S$ belongs to the finite strain space $\mathcal{B}$ which is defined as follows. Given a vector field $u\in
H^1(S,\mathbb{R}^3),$ by $\sym\nabla u$ we mean the bilinear form on $S_x,$ given by $\sym\nabla u(\tau,\eta)=\frac{1}{2}[\langle\partial_{\tau}u(x),\eta\rangle+\langle\partial_{\eta}u(x),\tau\rangle]$ for all $\tau,\eta\in S_x.$
Then the finite strain space is given by
\begin{equation*}
\mathcal{B}=\overline{\{\sym\nabla u^h:\quad u^h\in H^1(S,\mathbb{R}^3)\}}^{L^2(S)}
\end{equation*}
with the $L^2$ norm.

Next, we consider   the time-dependent case. Let the external force be given by
\be\label{eq1.15}
f^h(\xi,x+s\mathbf{n})=h\sqrt{e^h}f(h\xi,x),\quad f(\xi,x)\in L^2((0,\infty);L^2(S,\mathbb{R}^3)).
\ee
We assume that the initial data $\bar{w}^h$ and $\hat{w}^h$ have the following scaling conditions in terms of $h$
\begin{equation}\label{eq1.16}
\frac{1}{2}\int_{S^h}\mid\hat{w}^h(z)\mid^2dz+\int_{S^h}W(\nabla\bar{w}^h(z))dz\leqslant Che^h,
\end{equation}
where $C>0$ is a unform constant independent of $h.$

Let $u^h$ be a sequence of solutions to (\ref{eq1.5}) on $[0,T/h]\times S^h.$ As usual, we rescale $S^h$ to the fixed domain $S^{h_0}$ and  the time to  $t=h\xi.$
We set
\begin{equation}
y^h(t,x+s\mathbf{n}(x))\triangleq u^h(\frac{t}{h},x+\frac{sh}{h_0}\mathbf{n}(x)),\quad\mbox{on}\quad (0,T)\times S^{h_0}.\label{xxx1.19}
\end{equation}
It follows from (\ref{eq1.6}) that
\be\label{eq1.17}
y^h(0,x+s\mathbf{n}(x))=\bar{w}^h(x+\frac{sh}{h_0}\mathbf{n}(x)),\quad \partial_ty^h(0,x+s\mathbf{n}(x))=\frac{1}{h}\hat{w}^h(x+\frac{sh}{h_0}\mathbf{n}(x)),
\ee for $x+s\n\in S^{h_0}.$

We have the following.

\begin{thm}\label{thm1.1}
Let the assumptions $(\ref{eq1.1})-(\ref{eq1.23})$  and $(\ref{eq1.11})$ hold. Let $(\hat{w}^h)
\subset L^2(S^h,\mathbb{R}^3)$ and $(\bar{w}^h)\subset H^1(S^h,\mathbb{R}^3)$ satisfying the  boundary conditions $(\ref{eq1.7})$ and $(\ref{eq1.8})$ be the sequences of initial data of $(\ref{eq1.5})$
with the scaling assumption $(\ref{eq1.16}).$  Let   $h_0>0$ be given small and $T>0.$  Let the external force $f^h$ have the property $(\ref{eq1.15}).$ For all $h\in(0,h_0),$ let $y^h\in L^2((0,T);H^1(S^{h_0},\mathbb{R}^3))$ with
\begin{equation}
\partial_ty^h\in L^2((0,T);L^2(S^{h_0},\mathbb{R}^3)),\quad \partial^2_ty^h\in L^2((0,T);H^{-1}(S^{h_0},\mathbb{R}^3))\label{x1.25}
\end{equation}
be weak solutions to $(\ref{eq1.5})$ in $(0,T)\times S^{h_0}$ with initial data $(\ref{eq1.17}),$ the  boundary conditions $(\ref{eq1.18}),(\ref{eq1.19}),$ and the energy inequalities
\beq\label{eq1.24}
&&\int_{S^{h_0}}[\frac{h^2}{2}\mid\partial_ty^h(t,x+s\mathbf{n}(x))\mid^2+W(\nabla_hy^h(t,x+s\mathbf{n}(x)))]dz\nonumber \\
&&\leqslant \int_{S^{h_0}}[\frac12\mid\hat{w}^h(x+\frac{sh}{h_0}\mathbf{n}(x))\mid^2+W(\nabla_h\bar{w}^h(x+\frac{sh}{h_0}\mathbf{n}(x)))]\frac{\det F(\frac{sh}{h_0})}{\det F(s)} dz\nonumber \\
&&+h\sqrt{e^h}\int_0^t\int_{S^{h_0}}\langle f(t,x),\partial_ty^h\rangle dtdz,
\eeq
for $t\in(0,T),$  where $F(s)$ is given by
\be F(s)=\id+s\Pi.\label{xx1.17}\ee

Then for $y^h(t,x+s\mathbf{n}(x)),$
defined on the common domain $(0,T)\times S^{h_0},$ we have:

$(i)$\,\,\, $y^h$ converges in $L^{\infty}_{}((0,T),H^1(S^{h_0},\mathbb{R}^3))$ to $\pi.$

$(ii)$\,\,\, The scaled average displacements:
\begin{equation*}
V^h(t,x)=\frac{h}{\sqrt{e^h}}\fint^{\frac{h_0}{2}}_{-\frac{h_0}{2}}y^h(t,x+s\mathbf{n}(x))-xds
\end{equation*}
converges (up to a subsequence) in $L^q_{}((0,T),H^1(S,\mathbb{R}^3))$ to some $V\in L^{\infty}_{}((0,T),\mathcal{V})$ for $1\leqslant q<\infty.$ Besides,
$\pl_tV^h$ converges weakly-star in $L^{\infty}_{}((0,T),L^2(S,\mathbb{R}^3))$ to $\pl_tV$ and $V^h$ converges to $V$ in $L^{\infty}_{}((0,T),L^2(S,\mathbb{R}^3)),$
 respectively. Then $V\in W^{1,\infty}_{}((0,T),L^2(S,\mathbb{R}^3))\cap L^{\infty}_{}((0,T),\mathcal{V}).$

$(iii)$\,\,\, $\dfrac{1}{h}\sym\nabla V^h$ converges weakly in $L^2_{}((0,T),L^2(S))$ to some $B_{\tan}\in L^2_{}((0,T), \mathcal{B}).$

$(iv)$\,\,\, The couple $(V,B_{\tan})$ satisfies the following two variational dynamical equations. If $\kappa>0,$ for all $\tilde{V}\in L^2((0,T);\mathcal{V}\cap H^2_0(S,\mathbb{R}^3))\cap H^1_0((0,T);H^1_0(S,\mathbb{R}^3))$
with $\tilde{A}=\nabla \tilde{V}$ given as in $(\ref{eq1.12})$ and all $\tilde{B}_{\tan}\in  L^2_{}((0,T), \mathcal{B})$ there hold:
\begin{equation}\label{eq1.25}
\int_0^{T}\int_S\mathcal{L}_2(x,(B-\frac{\sqrt{\kappa}}{2}A^2)_{\tan}):\tilde{B}_{\tan}dxdt=0,
\end{equation}
\begin{eqnarray}\label{eq1.26}
&&\int_0^{T}\int_{S}\langle f,\tilde{V}\rangle dxdt+\int_0^{T}\int_{S}\langle V_t,\tilde{V}_t\rangle dxdt
=-\sqrt{\kappa}\int_0^{T}\int_{S}\mathcal{L}_2(x,(B-\frac{\sqrt{\kappa}}{2}A^2)_{\tan}):(A\tilde{A})_{\tan}dxdt \nonumber \\
&&+\frac{1}{12}\int_0^{T}\int_{S}\mathcal{L}_2(x,(\nabla(A\mathbf{n})-A\Pi)_{\tan}):[(\nabla(\tilde{A}\mathbf{n}))_{\tan}-(\tilde{A}\Pi)_{\tan}]dxdt.
\end{eqnarray}
If $\kappa=0,$ then  $(\ref{eq1.26})$ is still true where the first term in the right hand side of $(\ref{eq1.26})$    equals the zero.  Moreover, the initial data $V(0,x)=\bar{w}(x)\in\mathcal{V}$ with $\bar{w}(x)=0$ on $\partial S$ and $V_t(0,x)=\hat{w}(x)\in L^2(S)$ in the both cases, where $\bar{w}(x)$ and $\hat{w}(x)$ are the limits of $V^h(0,x)$ and $V^h_t(0,x)$ in a certain sense, respectively. The boundary values of $V$ satisfy
that $V(t,x)=0$ and $(\nabla V(t,x))^T\mathbf{n}=0$ for a.e. $(t,x)\in (0,T)\times\partial S.$
\end{thm}

\begin{rem} In Theorem $\ref{thm1.1}$ we have made the regularity assumption $(\ref{x1.25}).$ In the case of the thin plates such regularities have been established in $\cite{AMM1}.$
\end{rem}

\begin{rem}\label{rem1.1}
 By the scaling conditions $(\ref{eq1.16})$ for the initial data and from $\cite{LMP2},$ we have that
\begin{equation}\label{eq1.27}
\frac{1}{\sqrt{e^h}}\fint^{\frac{h_0}{2}}_{-\frac{h_0}{2}}\hat{w}^h(x+\frac{sh}{h_0}\mathbf{n})ds\rightharpoonup\hat{w}(x)\quad\mbox{in}
\quad L^2(S),
\end{equation}
and
\begin{equation}\label{eq1.28}
\frac{h}{\sqrt{e^h}}\fint^{\frac{h_0}{2}}_{-\frac{h_0}{2}}\bar{w}^h(x+\frac{sh}{h_0}\mathbf{n})-xds\rightarrow\bar{w}(x)
\quad\mbox{in}\quad H^1(S).
\end{equation}
Moreover, we also obtain that $\bar{w}(x)\in\mathcal{V}$ with $\bar{w}(x)=0$ on $\partial S.$ For more detail, see $\cite{AMM2,LMP2}.$
\end{rem}

\section{Some Modifications in the Stationary Shell Theory}
\setcounter{equation}{0}
\hskip\parindent
 We here list the some results in \cite{LMP2}.

\begin{thm}$\cite{LMP2}$\label{thm2.1}
Let $u^h\in H^1(S^h,\mathbb{R}^3)$ be a sequence of deformations of the thin shell $S^h.$ Assume $(\ref{eq1.11})$ and let the scaled energy $\frac{I^h(u^h)}{e^h}$ be uniformly bounded. Then there exists a sequence of matrix fields $R^h\in H^1(S,\mathbb{R}^3)$ with $R^h(x)\in \SO(3)$ for a.e. $x\in S,$ such that:
\begin{equation*}
\parallel\nabla u^h-R^h\pi\parallel_{L^2(S^h)}\leqslant Ch^{\frac{1}{2}}\sqrt{e^h}\quad\mbox{and}\quad \parallel\nabla R^h\parallel_{L^2(S)}\leqslant Ch^{-1}\sqrt{e^h}
\end{equation*}
and another sequence of matrices $Q^h\in \SO(3)$ such that

$(i)$ $\parallel(Q^h)^TR^h-\id\parallel_{L^p(S)}\leqslant C\frac{\sqrt{e^h}}{h},\quad\mbox{for}\quad p\in[1,\infty);$

$(ii)$ $\frac{h}{\sqrt{e^h}}((Q^h)^TR^h-\id)$ converges (up to a subsequence) to a skew-symmetric matrix field $\tilde{A},$ weakly in $H^1(S)$ and strongly in $L^p(S),$ where $p\in[1,\infty).$

Moreover, there is a sequence $c^h\in\mathbb{R}^3$ such that for the normalized rescaled deformations:
\begin{equation*}
\tilde{y}^h(x+s\mathbf{n})=(Q^h)^Ty^h(x+s\mathbf{n})-c^h,\quad\mbox{where}\quad y^h(x+s\mathbf{n})\triangleq u^h(x+\frac{sh}{h_0}\mathbf{n}(x))
\end{equation*}
defined on the common domain $S^{h_0},$ the following holds:

$(iii)$ $\parallel\nabla_hy^h-R^h\pi\parallel_{L^2(S^{h_0})}\leqslant C\sqrt{e^h}$ and $\tilde{y}^h$ converge in $H^1(S^{h_0})$ to $\pi;$

$(iv)$ The scaled average displacements $\tilde{V}^h,$ defined as $\tilde{V}^h(x)=\frac{h}{\sqrt{e^h}}\fint^{\frac{h_0}{2}}_{-\frac{h_0}{2}}\tilde{y}^h(x+s\mathbf{n}(x))-xds$ converge (up to a subsequence) in $H^1(S)$ to some $\tilde{V}\in\mathcal{V},$ whose gradient is given by $\tilde{A},$ as in $(\ref{eq1.11})$ and
\begin{equation*}
\lim_{h\rightarrow 0}\frac{h}{\sqrt{e^h}}((Q^h)^T\nabla_hy^h-\id)=\tilde{A}\pi,
\end{equation*}
in $L^2(S^{h_0})$ up to a subsequence;

$(v)$ $\frac{1}{h}sym\nabla \tilde{V}^h$ converges (up to a subsequence) in $L^2(S)$ to some symmetric matrix field $\tilde{B}_{\tan}\in\mathcal{B};$

$(vi)$
\begin{equation*}
\lim_{h\rightarrow 0}\frac{h^2}{e^h}\sym((Q^h)^TR^h-\id)=\frac{1}{2}\tilde{A}^2,\quad\mbox{in}\quad L^p(S),\quad\mbox{where}\quad p\in[1,\infty);
\end{equation*}

$(vii)$ Let $G^h=\frac{1}{\sqrt{e^h}}((R^h)^T\nabla_hy^h-\id).$ Then $G^h$ has a subsequence converging weakly in $L^2(S^{h_0})$ to a matrix field $G.$ Further,
the tangential minor of G satisfies that
\begin{equation*}
G(x+s\mathbf{n})\tau=G_0(x)\tau+\frac{s}{h_0}(\nabla(\tilde{A}\mathbf{n})-\tilde{A}\Pi)\tau,\quad \forall \tau\in S_x,
\end{equation*}
where $G_0(x)=\fint_{-\frac{h_0}{2}}^{\frac{h_0}{2}}G(x+s\mathbf{n})ds.$
\end{thm}

We observe that there is a byproduct $Q^h$ in Theorem \ref{thm2.1}. The construction of $Q^h$ in the appendices in \cite{LMP2} implies that it only depends on $h$ in the stationary case while in the time dependent case, it may depend on the time $t$ and be not differentiable on $t,$ which makes it more complicated in our analysis. In order to cope with it, we eliminate $Q^h$ by some idea  in \cite[Lemma 13]{LeM}.

We define the first moment by
\begin{equation*}
\tilde{\zeta}^h(x)=\fint_{-\frac{h_0}{2}}^{\frac{h_0}{2}}s[\tilde{y}^h(x+s\mathbf{n})-(x+\frac{sh}{h_0}\mathbf{n})]ds
\end{equation*}
 to determine the limit  of $\frac{1}{\sqrt{e^h}}\tilde{\zeta}^h$ as $h\rightarrow 0,$  which is useful for dealing with the related boundary value problem.

\begin{pro}\label{cor2.2}
Under the assumptions of Theorem $\ref{thm2.1},$ we have
\begin{equation}\label{eq2.1}
\frac{1}{\sqrt{e^h}}\tilde{\zeta}^h\rightharpoonup\frac{h_0}{12}\tilde{A}\mathbf{n},\quad\mbox{in}\quad H^1(S,\mathbb{R}^3),\quad\mbox{as}\quad h\rightarrow 0.
\end{equation}
\end{pro}
{\bf Proof} \,\,\, As in \cite{FGM2}, we set
\begin{equation*}
Y^h=\tilde{y}^h(x+s\mathbf{n})-(x+\frac{sh}{h_0}\mathbf{n}),\quad \bar{Y}^h=\fint_{-\frac{h_0}{2}}^{\frac{h_0}{2}}Y^hds,\quad Z^h=Y^h-\bar{Y}^h.
\end{equation*}
Thus, we have
$$\frac{h_0}{h}\partial_{\mathbf{n}}Z^h=\nabla_h\tilde{y}^h(x+s\mathbf{n}(x))\mathbf{n}(x)-\mathbf{n}(x).$$
Therefore, by $(iii)$ in Theorem \ref{thm2.1}, we obtain
\begin{equation*}
\parallel\frac{h_0}{h}\partial_{\mathbf{n}}Z^h-[(Q^h)^TR^h-\id]\mathbf{n}\parallel_{L^2(S^{h_0})}\leqslant C\sqrt{e^h}.
\end{equation*}

Since $\fint_{-\frac{h_0}{2}}^{\frac{h_0}{2}}Z^h=0$ and $\fint_{-\frac{h_0}{2}}^{\frac{h_0}{2}}s[(Q^h)^TR^h-\id]\mathbf{n}ds=0,$ by Poinc\'{a}re's inequality,
\begin{equation*}
\parallel\frac{h_0}{h}Z^h-s[(Q^h)^TR^h-\id]\mathbf{n}\parallel_{L^2(S^{h_0})}\leqslant C\sqrt{e^h}.
\end{equation*}

Multiply the quantity inside the above norm by $\frac{hs}{\sqrt{e^h}}$ and integrate with respect to $s$ over $(-h_0/2,h_0/2)$ to lead to
\begin{equation*}
\parallel\frac{1}{\sqrt{e^h}}\tilde{\zeta}^h-\frac{h_0}{12}\frac{h}{\sqrt{e^h}}[(Q^h)^TR^h-\id]\mathbf{n}\parallel_{L^2(S)}\leqslant Ch.
\end{equation*}
Thus, we have
$$\frac{1}{\sqrt{e^h}}\tilde{\zeta}^h\rightarrow\frac{h_0}{12}\tilde{A}\mathbf{n}\quad\mbox{in}\quad L^2(S).$$

Moreover, by straightforward calculation, we have for any $\tau\in S_x,$
\begin{equation*}
\partial_{\tau}\tilde{\zeta}^h(x)=\fint_{-\frac{h_0}{2}}^{\frac{h_0}{2}}s[\nabla_h\tilde{y}^h(x+s\mathbf{n}(x))-\id]F(\frac{sh}{h_0})\tau ds.
\end{equation*}

Using $(ii)$ and $(iii)$ in Theorem \ref{thm2.1}, we conclude that
$\frac{\tilde{\zeta}^h}{\sqrt{e^h}}$ are bounded in $H^1(S,\mathbb{R}^3).$ The proof is complete. \hfill$\Box$\\

Let
\begin{equation*}
V^h(x)=V^h[y^h](x)=\frac{h}{\sqrt{e^h}}\fint^{\frac{h_0}{2}}_{-\frac{h_0}{2}}y^h(x+s\mathbf{n}(x))-xds
\end{equation*}
and
\begin{equation*}
\zeta^h(x)=\fint^{\frac{h_0}{2}}_{-\frac{h_0}{2}}s[y^h(x+s\mathbf{n}(x))-(x+\frac{hs}{h_0}\mathbf{n}(x))]ds.
\end{equation*}
Now we consider some assumptions on the boundary value of $V^h$ and $\frac{1}{\sqrt{e^h}}\zeta^h$ and by these boundary value conditions, we will have that $Q^h=\id$ and $c^h=0.$

\begin{lem}\label{lem2.1} Let the assumptions in Theorem $\ref{thm2.1}$ hold and let the boundary conditions in $(\ref{eq1.18})$ be true.
Then the following asymptotic identities
\begin{equation}
Q^h=\id+O(\frac{\sqrt{e^h}}{h}),\quad \sym(Q^h-\id)=O(\frac{e^h}{h^2}),\quad c^h=O(\frac{\sqrt{e^h}}{h}). \label{x2.6}
\end{equation} hold.
\end{lem}

{\bf Proof}\,\,\,
First, we shall show that there is an open segment  $\Ga$ of $\pl S$ such that, for $x\in\Ga,$
\be x=x_{\tan}+\<x,\mathbf{n}(x)\>\n(x),\quad \int_\Ga xd\Ga=0,\quad \int_\Ga|x_{\tan}|d\Ga>0.\label{x2.5}\ee
In fact, we may assume that
$$\int_{\pl S}xdx=0.$$ Otherwise, we can translate $S.$
If $\int_{\pl S}|x_{\tan}|d\Ga>0,$ then we can let $\Ga=\pl S.$ Let
$$x_{\tan}=0\qfq x\in\pl S.$$
Then $\pl S$ is a curve on a sphere centered at the origin. Let $\Ga\subset\pl S$ be a segment such that
$$\int_\Ga xd\Ga\not=0.$$ We translate $S$ such that
$$\int_\Ga xd\Ga=0.$$ Since $\Ga$ is on a sphere not centered at the origin, we have
$$\int_\Ga |x_{\tan}|d\Ga>0.$$

Then
Comparing  the definitions of $V^h$ with $\tilde{V}^h$ and $\zeta^h$ with $\tilde{\zeta}^h,$ respectively, we obtain that
\begin{eqnarray}
&&\frac{\sqrt{e^h}}{h}V^h=\frac{\sqrt{e^h}}{h}Q^h\tilde{V}^h+(Q^h-\id)x+c^h,\label{eq2.5} \\
&&\zeta^h=Q^h\tilde{\zeta}^h+\frac{h_0}{12}h(Q^h-\id)\mathbf{n}(x), \label{eq2.6}
\end{eqnarray}
where we still denote $Q^hc^h$ by $c^h.$

Using (\ref{eq2.1}), (\ref{eq2.6}), the embedding $H^1(S)\hookrightarrow L^2(\Gamma),$ and  $\zeta^h=0$ in $L^2(\Gamma),$ we see that
\begin{equation}\label{eq2.7}
\parallel(Q^h-\id)\mathbf{n}\parallel_{L^2(\Gamma)}\leqslant C\frac{\sqrt{e^h}}{h},
\end{equation}
which yield by $Q^h\in \SO(3)$  that
\begin{equation}\label{eq2.8}
\parallel[(Q^h)^T-\id]\mathbf{n}\parallel_{L^2(\Gamma)}\leqslant C\frac{\sqrt{e^h}}{h}.
\end{equation}

We fix a point $x_0\in\Ga$ such that
 $$x_{0\tan}\not=0.$$
 Let $\tau_1(x),$ $\tau_2(x),$ $\mathbf{n}(x)$ be a local frame at $x_0$ on $\overline{S}$ with the positive orientation, where $\mathbf{n}(x)=\tau_1(x)\wedge\tau_2(x).$
Let $\Ga_0\subset\Ga$ be an open neighborhood of $x_0$ in $\Ga$ such that the frame is well defined on $\Ga_0.$ Let
$$Q_0(x)=\Big(\tau_1(x),\tau_2(x),\mathbf{n}(x)\Big)\qfq x\in\Ga_0.$$Obviously, we have $Q_0\in \SO(3).$ Let $Q^h_{\tan}$ denote the $2\times 2$ submatrix $(\langle Q^h\tau_i,\tau_j\rangle)_{i,j=1,2}$ of $Q_0^TQ^hQ_0.$
Via (\ref{eq2.7}) and (\ref{eq2.8}),  there exists a matrix $\hat{Q}^h(x)\in SO(2)$ such that
\begin{equation}\label{eq2.9}
\mid Q^h_{\tan}(x)-\hat{Q}^h(x)\mid\leqslant C\frac{\sqrt{e^h}}{h}\qfq x\in\Ga_0,
\end{equation}
where the constant $C$ is independent of $x\in\Gamma.$

From (\ref{eq2.5}), (\ref{eq2.7}), and  $V^h=0$ on $\Ga$ it follows that
\begin{equation}\label{eq2.10}
\parallel(Q^h-\id)x+c^h\parallel_{L^2(\Gamma)}\leqslant C\frac{\sqrt{e^h}}{h}.
\end{equation}
It follows  from (\ref{x2.5}), (\ref{eq2.8}), and (\ref{eq2.10}) that
\begin{eqnarray}
&&\mid c^h\mid\leqslant C\frac{\sqrt{e^h}}{h},\quad\parallel(Q^h-\id)x_{\tan}\parallel_{L^2(\Gamma)}\leqslant C\frac{\sqrt{e^h}}{h}. \label{eq2.12}
\end{eqnarray}
By  (\ref{eq2.9}) and (\ref{eq2.12}), we have
\begin{equation}\label{eq2.13}
\parallel(\hat{Q}^h(x)-\id)x_{\tan}\parallel_{L^2(\Gamma_0)}\leqslant C\frac{\sqrt{e^h}}{h}.
\end{equation}
Now, since $\hat Q^h(x)\in SO(2),$ $2|\hat Q^h-\id)x_{\tan}|^2=|\hat Q^h(x)-\id|^2|x_{\tan}|^2.$ It follows from  (\ref{eq2.13}) that
\begin{equation}\label{eq2.14}
\parallel Q^h_{\tan}-\id\parallel_{L^2(\Gamma_0)}\leqslant C\frac{\sqrt{e^h}}{h}.
\end{equation}

From (\ref{eq2.9})  and (\ref{eq2.14}), we obtain
\begin{equation}\label{eq2.15}
\mid Q^h-\id\mid\leqslant C\frac{\sqrt{e^h}}{h}.
\end{equation}
Moreover, from the relation
\begin{equation*}
2\sym(Q^h-\id)=-((Q^h)^T-\id)(Q^h-\id),
\end{equation*} we have the second asymptotic identity in (\ref{x2.6}).  \hfill$\Box$\\

It follows from (\ref{x2.6}), (\ref{eq2.5}), and (\ref{eq2.6}) that
\begin{lem}\label{lem2.2}
Suppose that
\begin{equation}\label{eq2.2}
I^h(y^h)=\int_{S^{h_0}}W(\nabla_hy^h)dz\leqslant Ce^h,\quad\lim_{h\rightarrow 0}\frac{e^h}{h^2}=0.
\end{equation}
Moreover, let $(\ref{eq1.18})$ and $(\ref{x2.5})$ hold. Then
\begin{eqnarray}
&&V^h\rightarrow V\quad\mbox{in}\quad H^1(S,\mathbb{R}^3)\quad\mbox{with}\quad V\in\mathcal{V},\label{eq2.3} \\
&&\frac{1}{\sqrt{e^h}}\zeta^h\rightharpoonup \frac{h_0}{12}A\mathbf{n}\quad\mbox{in}\quad H^1(S,\mathbb{R}^3).\label{eq2.4}
\end{eqnarray} Moreover, there is $A_0\in\SO(3)$ such that
$$V=\tilde{V}+A_0x+c_0,\quad A=\tilde{A}+A_0\in \so(3),$$ where $\so(3)$ is the set of all the $3\times3$ anti-symmetric matrices.
\end{lem}

\begin{rem} In Lemma $\ref{lem2.1},$
$$A_0=\lim_{h\rightarrow0}\frac{h}{\sqrt{e^h}}(Q^h-\id),\quad c_0=\lim_{h\rightarrow0}\frac{h}{\sqrt{e^h}}c^h.$$
\end{rem}

Now, by applying Lemmas \ref{lem2.1} and \ref{lem2.2}, Theorem \ref{thm2.1} can be rewritten as the following.

\begin{thm}$\cite{LMP2}$\label{thm2.2}
Let $u^h\in H^1(S^h,\mathbb{R}^3)$ be a sequence of deformations of thin shell $S^h.$ Suppose that $(\ref{eq1.11})$ and all the assumptions in Lemma $\ref{lem2.1}$ hold true. Then there is a sequence of matrix fields $R^h\in H^1(S,\mathbb{R}^3)$ with $R^h(x)\in \SO(3)$ for a.e. $x\in S,$ satisfying:
\be\parallel\nabla u^h-R^h\pi\parallel_{L^2(S^h)}\leqslant Ch^{\frac{1}{2}}\sqrt{e^h}, \quad\parallel\nabla R^h\parallel_{L^2(S)}\leqslant Ch^{-1}\sqrt{e^h};
\ee

$(i)$\,\,\, $\parallel R^h-\id\parallel_{H^1(S)}\leqslant C\frac{\sqrt{e^h}}{h};$

$(ii)$\,\,\, $A^h\triangleq\dfrac{h}{\sqrt{e^h}}(R^h-\id)$ converges (up to a subsequence) to a skew-symmetric matrix field $A=\tilde{A}+A_0,$ weakly in $H^1(S)$ and strongly in $L^p(S).$

Moreover, for the rescaled deformations
\begin{equation*}
y^h(x+s\mathbf{n})\triangleq u^h(x+\frac{sh}{h_0}\mathbf{n}(x))
\end{equation*}
defined on the common domain $S^{h_0},$ the following holds:

$(iii)$ $\parallel\nabla_hy^h-R^h\pi\parallel_{L^2(S^{h_0})}\leqslant C\sqrt{e^h}$ and $y^h$ converges in $H^1(S^{h_0})$ to $\pi;$

$(iv)$ The scaled average displacements $V^h,$ defined as $V^h(x)=\frac{h}{\sqrt{e^h}}\fint^{\frac{h_0}{2}}_{-\frac{h_0}{2}}y^h(x+s\mathbf{n}(x))-xds$ converges (up to a subsequence) in $H^1(S)$ to $V=\tilde{V}+A_0x+c_0\in\mathcal{V},$  and
\begin{equation}\label{eq2.19}
\lim_{h\rightarrow 0}\dfrac{h}{\sqrt{e^h}}(\nabla_hy^h-\id)=A\pi\quad\mbox{in}\quad L^2(S^{h_0});
\end{equation}

$(v)$\,\,\, $\dfrac{1}{h}\sym\nabla V^h$ converges (up to a subsequence) in $L^2(S)$ to some $B_{\tan}\in\mathcal{B};$

$(vi)$
\begin{equation*}
\lim_{h\rightarrow 0}\frac{h^2}{e^h}\sym(R^h-\id)=\frac{1}{2}A^2\quad\mbox{in}\quad L^p(S),\quad\mbox{where}\quad p\in[1,\infty);
\end{equation*}

$(vii)$\,\,\, Let $G^h=\dfrac{1}{\sqrt{e^h}}((R^h)^T\nabla_hy^h-\id).$ Then $G^h$ has a subsequence converging weakly in $L^2(S^{h_0})$ to a matrix field $G.$ Further,
\begin{equation*}
G(x+s\mathbf{n})\tau=G_0(x)\tau+\frac{t}{h_0}(\nabla(A\mathbf{n})-A\Pi)\tau,\quad \forall \tau\in S_x,
\end{equation*}
where $G_0(x)=\fint_{-\frac{h_0}{2}}^{\frac{h_0}{2}}G(x+s\mathbf{n})ds.$
\end{thm}

\section{The Proof of Theorem \ref{thm1.1}}
\setcounter{equation}{0}
\hskip\parindent
We need to make some preparations for  deriving the two-dimensional evolutionary nonlinear shell model from the corresponding three-dimensional elastodynamic system.

Let $y^h$ be given in (\ref{xxx1.19}). Define
\begin{equation*}
\nabla_hy^h(t,x+s\mathbf{n}(x))=\nabla u^h(\frac{t}{h},x+\frac{sh}{h_0}\mathbf{n}(x)).
\end{equation*}
A straightforward calculation yields, for all $x\in S,$ $s\in(-\frac{h_0}{2},\frac{h_0}{2}),$ and $\tau\in S_x,$
\beq
&&\partial_{\tau}y^h(t,x+s\mathbf{n}(x))=\nabla_hy^h(t,x+s\mathbf{n}(x))F(\frac{sh}{h_0})F^{-1}(s)\tau,\quad\quad \quad\label{x1.19}\\
&&\partial_{\mathbf{n}}y^h(t,x+s\mathbf{n}(x))=\frac{h}{h_0}\nabla_hy^h(t,x+s\mathbf{n}(x))\mathbf{n}(x),\nonumber
\eeq where $F(s)$ is given in (\ref{xx1.17}).
Therefore, the boundary conditions in (\ref{eq1.7}) and (\ref{eq1.8}) become
\begin{eqnarray}
&&y^h(t,x+s\mathbf{n}(x))=x+\frac{sh}{h_0}\mathbf{n}(x)\quad\mbox{on}\quad \{x+s\mathbf{n}(x): x\in\partial S, s\in(-\frac{h_0}{2},\frac{h_0}{2})\}\times(0,T),\label{eq1.18} \\
&&DW(\nabla_hy^h)\mathbf{n}=0\quad\mbox{on}\quad \{x\pm\frac{h_0}{2}\mathbf{n}(x): x\in S\}\times(0,T).\label{eq1.19}
\end{eqnarray}

The conditions (\ref{eq1.16}) are
\beq\label{eq1.20}
\int_{S^{h_0}}[\frac{h^2}{2}|\pl_ty^h(0,x+s\mathbf{n}(x))\mid^2+W(\nabla_hy^h(0,x+s\mathbf{n}(x)))]\frac{\det F(\frac{sh}{h_0})}{\det F(s)}dz\leqslant Ch_0e^h.
\eeq

For each $\varphi\in C_0^\infty((0,\infty)\times S^{h_0}),$ we consider the test function
$$\psi^h(\xi,x+s\n)=\varphi(h\xi,x+\frac{sh_0}{h}\n)\qfq (\xi,x+s\n)\in(0,\infty)\times S^h.$$
We have the following Euler-Lagrange equations
\be\label{eq1.21}
\int_0^{T/h}\int_{S^{h}}[\langle u^h_\xi,\psi^h_\xi\rangle-DW(\nabla u^h):\nabla\psi^h+\langle f^h,\psi^h\rangle]d\xi dz=0.
\ee
Let $\tau_1,$ $\tau_2$ be a local form on $S.$  In (\ref{eq1.21}) $\nabla\psi^h$ is given by
$$\nabla\psi^h(\xi,x+s\n)\tau_i=\nabla\varphi(h\xi,x+\frac{sh_0}{h}\n)F(\frac{sh}{h_0})F^{-1}(s)\tau_i\qfq i=1,\,2,$$
$$\nabla_{\n}\psi^h(\xi,x+s\n)=\frac{h_0}{h}\nabla_{\n}\varphi(h\xi,x+\frac{sh_0}{h}\n),$$ where $F(s)$ is given in (\ref{xx1.17}).
It is easy to check that (\ref{eq1.21}) can be rewritten as
\begin{equation}\label{eq1.22}
\int_{T,S,h_0}[\langle hy^h_t,h\varphi_t\rangle-DW(\nabla_hy^h):\nabla\varphi P_h+h\sqrt{e^h}\langle f,\varphi\rangle]\det F(\frac{sh}{h_0})dsdxdt=0,
\end{equation} where
\be\int_{T,S,h_0}=\int_0^{T}\int_{S}\fint^{\frac{h_0}{2}}_{-\frac{h_0}{2}},\label{xn1.23}\ee
$$P_h=\Big(F^{-1}(\frac{sh}{h_0})F(s)\tau_1, F^{-1}(\frac{sh}{h_0})F(s)\tau_2, \frac{h_0}{h}\n\Big)\Big(\tau_1,\tau_2,\n\Big)^T.$$

By similar arguments as in \cite{AMM2} and \cite{LMP2}, we have Lemma \ref{l3.1} below.
\begin{lem}\label{l3.1} $(i),$ $(ii),$ and $(iii)$ in Theorem $\ref{thm1.1}$ hold.
\end{lem}

\begin{lem}\label{l3.2} $(i)$\,\,\,$V^h(t,x)=\dfrac{h}{\sqrt{e^h}}\fint^{\frac{h_0}{2}}_{-\frac{h_0}{2}}y^h(t,x+s\mathbf{n}(x))-xds$ converges   to $V\in L^\infty((0,T),\mathcal{V})$ weakly-star in  $L^\infty((0,T),H^1(S,\R^3)),$ where
$$V(t,x)=0,\quad (\nabla V)^T\n=0\qfq (t,x)\in(0,T)\times\pl S,$$
$$V(0,x)=\bar{w}(x),\quad V_t(0,x)=\hat{w}(x)\qfq x\in S.$$

 $(ii)$\,\,\,$A^h(t,x)\rw A\in L^\infty((0,T),H^1(S,\R^{3\times3}))$ weakly-star in $L^\infty((0,T),H^1(S,\R^{3\times3})),$
 $\sym A^h\rightarrow 0$ strongly in $L^{\infty}((0,T),L^p(S,\mathbb{R}^{3\times 3})),$ for $1\leqslant p<\infty,$ and
 $(A^h\tau)$ is compact in $L^{q}((0,T),L^p(S,\mathbb{R}^{3}))$ for all $1\leqslant q<\infty, 2\leqslant p<\infty$ and $\tau\in \mathcal{X}(S),$
  where
$$A^h(t,x)=\frac{h}{\sqrt{e^h}}(R^h(t,x)-\id),\quad A_{\tan}=(\nabla V)_\tan,$$ and
$$\frac{h}{\sqrt{e^h}}(\sym A^h)_{\tan}\rightarrow\frac{1}{2}A^2_{\tan}\quad\mbox{strongly in}\quad L^{2}((0,T),L^2(S)).$$
Moreover, $A$ is a skew-symmetric matrix and the above $A$ satisfies
$$A(t,x)\n=0\qfq (t,x)\in(0,T)\times\pl S.$$

$(iii)$\,\,\,$G^h=\dfrac{1}{\sqrt{e^h}}((R^h)^T\nabla_hy^h-\id)\rightarrow G\quad\mbox{weakly-star in}\quad L^{\infty}((0,T),L^2(S^{h_0},\mathbb{R}^{3\times 3})),$ where
\be\sym G(t, x+s\mathbf{n}(x))_{\tan}=(B-\frac{\sqrt{\kappa}}{2}A^2)_{\tan}+\frac{s}{h_0}(\nabla(A\mathbf{n})-A\Pi)_{\tan},\label{xxxx3.7}\ee
$$B_{\tan}=\lim_{h\rw0}\frac{1}{h}\sym\nabla V^h\quad\mbox{weakly in}\quad L^2((0,T),L^2(S)).$$
\end{lem}

{\bf Proof}\,\,\,By  a similar argument as in \cite{AMM2} and \cite{LMP2}, we obtain (i)-(iii), where
$$G(t,x+s\mathbf{n})\tau=G_0(t,x)\tau+\frac{t}{h_0}(\nabla(A\mathbf{n})-A\Pi)\tau,\quad \forall \tau\in S_x,\quad G_0(t,x)=\fint_{-\frac{h_0}{2}}^{\frac{h_0}{2}}G(t,x+s\mathbf{n})ds.$$

Next, we compute $G_0(t,x).$ We have
\beq \frac{1}{h}\sym\nabla_{\tan} V^h(t,x)&&=\dfrac{1}{\sqrt{e^h}}\fint^{\frac{h_0}{2}}_{-\frac{h_0}{2}}\sym[\nabla_\tan y^h(t,x+s\mathbf{n}(x))F(s)-\id]ds\nonumber\\
&&=\dfrac{1}{\sqrt{e^h}}\fint^{\frac{h_0}{2}}_{-\frac{h_0}{2}}\sym[\nabla_\tan y^h(t,x+s\mathbf{n}(x))F(s)-(R^h)_\tan]ds\nonumber\\
&&\quad+\dfrac{1}{\sqrt{e^h}}\sym(R^h-\id)_\tan.\label{xxxx3.8}\eeq
It follows from (ii) that
$$\dfrac{1}{\sqrt{e^h}}\sym(R^h-\id)_\tan=\frac{\sqrt{e^h}}{h^2}\frac{h}{\sqrt{e^h}}\sym A_\tan^h\rw \frac{\sqrt{\kappa}}2A^2_\tan\quad\mbox{strongly in}\quad L^{2}((0,T),L^2(S)). $$
To treat the first term in the right hand side of (\ref{xxxx3.8}), we observe that
\beq &&\frac{1}{\sqrt{e^h}}[\nabla y^hF(s)-R^h]\tau=R^hG^h\tau+\frac{sh}{h_0\sqrt{e^h}}\nabla_hy^h\Pi\tau\qfq\tau\in S_x.\nonumber\eeq
Using the above formulas in (\ref{xxxx3.8}) and letting $h\rw0,$ we obtain
$$B_\tan=\sym[G_0(t,x)]_\tan+\frac{\sqrt{\kappa}}2A^2_\tan,$$ which yields the formula (\ref{xxxx3.7}). \hfill$\Box$

\begin{lem} \label{l3.3} Let
$$E^h=\frac{1}{\sqrt{e^h}}DW(\id+\sqrt{e^h}G^h).$$ Then

$(i)$\,\,\,$E^h, R^hE^h\rightarrow E=\mathcal{L}_3G$ weakly-star in $L^{\infty}((0,T),L^2(S^{h_0},\mathbb{R}^{3\times 3})),$ where $E$ is symmetric.

$(ii)$\,\,\,$E_{\tan}(t,x+s\mathbf{n}(x))=\mathcal{L}_2(x,G_{\tan}(t,x+s\mathbf{n}(x))).$

$(iii)$\,\,\,
$$\lim_{h\rightarrow 0}\frac{1}{h}\parallel \skew E^h\parallel_{L^{\infty}((0,T),L^p(S^{h_0}\mathbb{R}^{3\times 3}))}=0\quad\mbox{for}\quad p\in(1,4/3).$$

Moreover, let
\begin{equation*}
\bar{E}(t,x)=\fint^{\frac{h_0}{2}}_{-\frac{h_0}{2}}E(t,x+s\mathbf{n}(x))ds,\quad\hat{E}(t,x)=\fint^{\frac{h_0}{2}}_{-\frac{h_0}{2}}sE(t,x+s\mathbf{n}(x))ds.
\end{equation*}
Then

$(iv)$\,\,\,$\bar{E}_{\tan}(t,x)=\fint^{\frac{h_0}{2}}_{-\frac{h_0}{2}}\mathcal{L}_2(x,G_{\tan}(t,x+s\mathbf{n}(x)))ds=\mathcal{L}_2(x,(B-\frac{\sqrt{\kappa}}{2}A^2)_{\tan}).$

$(v)$\,\,\,$\hat{E}_{\tan}(t,x)=\fint^{\frac{h_0}{2}}_{-\frac{h_0}{2}}s\mathcal{L}_2(x,G_{\tan}(t,x+s\mathbf{n}(x)))ds=
\dfrac{h_0}{12}\mathcal{L}_2(x,(\nabla(A\mathbf{n})-A\Pi)_{\tan}).$
\end{lem}

{\bf Proof}\,\,\,(i)\,\,\,From (\ref{eq1.23}), $\{E^h\}$ is bounded in $L^{\infty}((0,T),L^2(S^{h_0},\mathbb{R}^{3\times 3})),$ and thus
$$E^h\rw E=\mathcal{L}_3G\quad\mbox{weakly-star in}\quad L^{\infty}((0,T),L^2(S^{h_0},\mathbb{R}^{3\times 3})),$$ arguing as in \cite[Proposition 2.2]{MP}.
By (i) in Theorem \ref{thm2.2} and the weakly-star convergence of $E^h,$ we also have
$$R^hE^h\rw E\quad\mbox{weakly-star in}\quad L^{\infty}((0,T),L^2(S^{h_0},\mathbb{R}^{3\times 3})),$$

(ii) follows from an argument as in  \cite[Lemma 2.3]{L}.

(iii)\,\,\,Since $DW(F)F^T$ is symmetric for all $F\in\mathbb{R}^{3\times 3}$(\cite[p.257]{AMM2}), we have
$$E^h-(E^h)^T+\sqrt{e^h}[E^h(G^h)^T-G^h(E^h)^T]=0.$$
It follows from (iii) in Lemma \ref{l3.2} and (i) that
$$\sup_{t\in(0,T)}\|\skew E^h\|_{L^2(S,\R^{3\times3})}\leq C,\quad\sup_{t\in(0,T)}\|\skew E^h\|_{L^1(S,\R^{3\times3})}\leq C\sqrt{e^h}.$$ By the interpolation inequality, we have for $p\in (1,2),$
\begin{equation*}
\frac{1}{h}\sup_{[0,T']}\parallel \skew E^h\parallel_{L^p(S^{h_0})}\leqslant\frac{1}{h}\sup_{[0,T']}\parallel \skew E^h\parallel_{L^1}^\theta\sup_{[0,T']}\parallel \skew E^h\parallel_{L^2}^{1-\theta}\leqslant\frac{C}{h}(e^h)^{\frac{\theta}{2}},
\end{equation*}
where $\frac{1}{p}=\theta+\frac{1-\theta}{2}$ and $\theta\in(1/2,1).$ Thus (iii) follows.

(iv) and (v) follow from (iii) in Lemma \ref{l3.2} and (ii), respectively.

Now, we are ready to prove Theorem \ref{thm1.1}.\\

{\bf Proof of Theorem \ref{thm1.1}} \,\,\, 1).\,\, {\bf Proof of (\ref{eq1.25})}\,\,\,Using the formulas $DW(F)=Q^TDW(QF)$ for $F\in\R^{3\times3},$ $Q\in SO(3),$ we have
$$ DW(\nabla_hy^h)=R^hDW(\id+\sqrt{e^h}G^h)=\sqrt{e^h}R^hE^h. $$ In addition, from (\ref{eq1.16}) and (\ref{eq1.24}), we have
\be\frac{h^2}{\sqrt{e^h}}\sup_{t\in[0,T]}\|y^h_t\|^2_{L^2(S^{h_0})}\leq C(1+\|f\|^2_{L^2((0,T)\times S)})\sqrt{e^h}.\label{xxx3.8}\ee

For any $\phi\in L^2((0,T),L^2(S^{h_0},\R^3)),$ let
$$\var(t,x+s\n))=\int_{-h_0/2}^s\phi(t,x+\eta\n)d\eta.$$ Then
$$\nabla_\n\var=\phi,\quad\nabla\var P_h\n=\frac{h_0}{h}\phi.$$ Using this $\var$ in (\ref{eq1.22}), we obtain
\beq h_0\int_{T,S,h_0}\<R^hE^h\n,\phi\>dsdxdt&&=\int_{T,S,h_0}[\langle\frac{h^2}{\sqrt{e^h}}y^h_t,h\varphi_t\rangle-h\sum_{i=1}^2\<R^hE^h\tau_i,\nabla\varphi P_h\tau_i\>\nonumber\\
&&\quad+h^2\langle f,\varphi\rangle]\det F(\frac{sh}{h_0})dsdxdt,\nonumber\eeq which yield, by letting $h\rw0,$
$$\int_{T,S,h_0}\<E\n,\phi\>dsdxdt=0\quad\mbox{for any}\quad \phi\in L^2((0,T),L^2(S^{h_0},\R^3)),$$ that is,
\be E\n=0\quad\mbox{a.e. on}\quad (0,T)\times S^{h_0}.\label{xxxx3.10}\ee

For $\phi(t,x)\in L^2((0,T);H^1(S,\mathbb{R}^3))\cap H_0^1((0,T);L^2(S,\mathbb{R}^3))$ with $\phi=0$ on $(0,T)\times\partial S,$ this time we let
$$\var(t,x+s\n)=\phi(t,x).$$ Then
$$ \nabla_\n\var=0,\quad\nabla\var P_h\tau_i=\nabla\phi F^{-1}(\frac{sh}{h_0})\tau_i,\quad i=1,\,2.$$
Thus (\ref{eq1.22}) can be written as
\beq\label{eq3.28}
\int_{T,S,h_0}[\langle \frac{h^2}{\sqrt{e^h}}y^h_t,\varphi_t\rangle
-\sum_{i=1}^2\<R^hE^h\tau_i:\nabla\phi F^{-1}(\frac{sh}{h_0})\tau_i\>+h\langle f,\varphi\rangle]\det F(\frac{sh}{h_0}) dsdxdt=0.\quad\quad\label{xxx3.7}
\eeq

From (i) in Theorem \ref{thm2.2}, (i) and (ii) in Lemma \ref{l3.3}, we obtain
\beq \lim_{h\rw0}\sum_{i=1}^2\<R^hE^h\tau_i:\nabla\phi F^{-1}(\frac{sh}{h_0})\tau_i\>&&=\sum_{i=1}^2\<E\tau_i,\nabla\phi\tau_i\>=E_{\tan}:\sym\nabla_{\tan}\phi\nonumber\\
&&=\mathcal{L}_2(x,G_{\tan}(t,x+s\mathbf{n}(x))):\sym\nabla_{\tan}\phi,\label{xxx3.9}\eeq where the convergence is  weakly-star in $L^{\infty}((0,T),L^2(S^{h_0},\mathbb{R})).$

Letting $h\rw0$ in (\ref{xxx3.7}) after using (\ref{xxx3.8}) and (\ref{xxx3.9}), we obtain (\ref{eq1.25}) by (iv) in Lemma \ref{l3.3}. \\

 2).\,\, {\bf Proof of (\ref{eq1.26})}\,\,\,
Let $\tilde{V}\in L^2((0,T),\mathcal{V}\cap H^2_0(S,\mathbb{R}^3))\cap H^1_0((0,T),H^1_0(S,\mathbb{R}^3))$ and let $\varphi(t,x+s\mathbf{n}(x))=s\tilde{A}(t,x)\mathbf{n}(x),$ where $\tilde{A}$ is skew-symmetric such that $\partial_{\tau}\tilde{V}=\tilde{A}\tau$ for all $\tau\in S_x.$ Then
$$\nabla_{\n}\var=\tilde{A}\n,\quad\nabla\var F^{-1}(\frac{sh}{h_0})F(s)\tau=s\nabla(\tilde{A}\n)F^{-1}(\frac{sh}{h})\tau\qfq\tau\in S_x.$$
Using (\ref{eq1.22}), (\ref{xxx3.8}), (i), (iv) and (v) in Lemma \ref{l3.3}, we obtain
\beq &&\lim_{h\rw0}\frac{h_0}{h}\int_{T,S,h_0}\<R^hE^h\n,\tilde{A}\n\>F(\frac{sh}{h_0})dsdxdt\nonumber\\
&&=\lim_{h\rw0}\int_{T,S,h_0}s[\langle \frac{h^2}{\sqrt{e^h}}y^h_t,\pl_t\tilde{A}\n\rangle-\sum_{i=1}^2\<R^hE^h\tau_i,\nabla(\tilde{A}\n)F^{-1}(\frac{sh}{h})\tau_i\>+h
\langle f,\tilde{A}\n\rangle]\det F(\frac{sh}{h_0})dsdxdt\nonumber\\
&&=-\int_{T,S,h_0}s\sum_{i=1}^2\<E\tau_i,\nabla(\tilde{A}\n)\tau_i\>dsdxdt=-\int_0^T\int_S\hat E_{\tan}:\nabla_{\tan}(\tilde{A}\n)dxdt\nonumber\\
&&=-\frac{h_0}{12}\int_0^T\int_S\mathcal{L}_2(x,(\nabla(A\mathbf{n})-A\Pi)_{\tan}):\nabla_{\tan}(\tilde{A}\n)dxdt.\label{xxx3.10}\eeq

Next, for $\tilde{V}(t,x)\in L^2((0,T),\mathcal{V}\cap H^2_0(S,\mathbb{R}^3))\cap H^1_0((0,T),H^1_0(S,\mathbb{R}^3)),$ let
$$\var(t,x+s\n)=\tilde{V}(t,x).$$ Let $\tilde{A}$ be the skew-symmetric matrix such that $\partial_{\tau}\tilde{V}=\tilde{A}\tau$ for all $\tau\in S_x.$ Then
$$\nabla_\n\var=\nabla\tilde{V}\n=0,\quad \nabla\var\tau=\tilde{A}F^{-1}(s)\tau\qfq\tau\in S_x.$$ It follows from (\ref{eq1.22}) that
\beq&&\int_{T,S,h_0}[\langle \frac{h}{\sqrt{e^h}}y^h_t,\tilde{V}_t\rangle+\langle f,\tilde{V}\rangle]\det F(\frac{sh}{h_0})dsdxdt\nonumber\\
&&=\frac{1}{h}\int_{T,S,h_0}\sum_{i=1}^2\<R^hE^h\tau_i,\tilde{A}F^{-1}(\frac{sh}{h_0})\tau_i\>\det F(\frac{sh}{h_0})dsdxdt\nonumber\\
&&=\frac{\sqrt{e^h}}{h^2}\int_{T,S,h_0}\sum_{i=1}^2\<A^hE^h\tau_i,\tilde{A}F^{-1}(\frac{sh}{h_0})\tau_i\>\det F(\frac{sh}{h_0})dsdxdt\nonumber\\
&&\quad+\frac{1}{h}\int_{T,S,h_0}\sum_{i=1}^2\<E^h\tau_i,\tilde{A}F^{-1}(\frac{sh}{h_0})\tau_i\>\det F(\frac{sh}{h_0})dsdxdt.\label{xxx3.14}\eeq
Let $x\in S$ be fixed. For simplicity, we select an orthonormal basis $\tau_1(x),$ $\tau_2(x)$ in $S_x$ such that
$$\Pi\tau_i=\nabla_{\tau_i}\n=\lam_i\tau_i\qfq i=1,\,2,$$ where $\lam_1\lam_2$ is the Gaussian curvature. Then
$$F^{-1}(\frac{sh}{h_0})\tau_i=\frac{h_0\tau_i}{h_0+sh\lam_i}\qfq i=1,\,2.$$
Thus we have
\beq
&&\sum_{i=1}^2\<E^h\tau_i,\tilde{A}F^{-1}(\frac{sh}{h_0})\tau_i\>=\sum_{i=1}^2\<E^h\tau_i,\tilde{A}\tau_i\>-h\sum_{i=1}^2\frac{\lam_i}{h_0+sh\lam_i}\<sE^h\tau_i,\tilde{A}\tau_i\>\nonumber\\
&&= E^h_\tan:\tilde{A}_\tan+\sum_{i=1}^2\<E^h\tau_i,\n\>\<\tilde{A}\tau_i,\n\>-h\sum_{i=1}^2\frac{\lam_i}{h_0+sh\lam_i}\<sE^h\tau_i,\tilde{A}\tau_i\>\nonumber\\
&&= \skew E^h_\tan:\tilde{A}_\tan-\<{E^h}^T\n,\tilde{A}\n\>-h\sum_{i=1}^2\frac{\lam_i}{h_0+sh\lam_i}\<sE^h\tau_i,\tilde{A}\tau_i\>\nonumber\\
&&= \skew E^h_\tan:\tilde{A}_\tan+2\<\skew E^h\n,\tilde{A}\n\>+\frac{\sqrt{e^h}}{h}\<A^hE^h\n,\tilde{A}\n\>-\<R^hE^h\n,\tilde{A}\n\>\nonumber\\
&&\quad-h\sum_{i=1}^2\frac{\lam_i}{h_0+sh\lam_i}\<sE^h\tau_i,\tilde{A}\tau_i\>,\nonumber\eeq since $\tilde{A}^T=-\tilde{A}$ and $\<\tilde{A}\n,\n\>=0.$
Using (iv) and (iii) in Lemma \ref{l3.3}, (\ref{xxx3.10}), and (\ref{eq1.11}), we obtain
\beq
&&\lim_{h\rw0}\frac{1}{h}\int_{T,S,h_0}\sum_{i=1}^2\<E^h\tau_i,\tilde{A}F^{-1}(\frac{sh}{h_0})\tau_i\>\det F(\frac{sh}{h_0})dsdxdt\nonumber\\
&&=\frac{1}{12}\int_0^T\int_S\mathcal{L}_2(x,(\nabla(A\mathbf{n})-A\Pi)_{\tan}):[\nabla(\tilde{A}\n)-\tilde{A}\Pi]_{\tan}dxdt.\label{xxx3.12}\eeq
Moreover, from (\ref{xxxx3.10}), (iv) and (i) in Lemma \ref{l3.3} and (ii) in Lemma \ref{l3.2} including the compactness of $(A^h\tau)$ and the strong convergence of $\sym A^h,$ we have
\beq&&\lim_{h\rw0}\int_{T,S,h_0}\sum_{i=1}^2\<A^hE^h\tau_i,\tilde{A}F^{-1}(\frac{sh}{h_0})\tau_i\>\det F(\frac{sh}{h_0})dsdxdt\nonumber\\
&&=-\int_0^T\int_S\sum_{i=1}^2\<\bar{E}\tau_i,A\tilde{A}\tau_i\>dxdt=-\int_0^T\int_S\bar{E}_{\tan}\:(A\tilde{A})_{\tan}dxdt\nonumber\\
&&=-\int_0^T\int_S\mathcal{L}_2(x,(B-\frac{\sqrt{\kappa}}{2}A^2)_{\tan}):(A\tilde{A})_{\tan}dxdt.\label{xxx3.16}\eeq

Finally, using (\ref{xxx3.8}), (\ref{xxx3.12}) and (\ref{xxx3.16}) in (\ref{xxx3.14}), we let $h\rw0$ to obtain (\ref{eq1.26}).  \hfill$\Box$

\end{document}